\documentclass[11pt,a4paper,reqno]{amsart}

\usepackage[T1]{fontenc}
\usepackage{lmodern}
\usepackage{amssymb,mathtools}
\usepackage{microtype}
\usepackage{graphicx}
\usepackage[margin=1in]{geometry}
\usepackage[hidelinks,hypertexnames=false]{hyperref}

\newtheorem{theorem}{Theorem}[section]
\newtheorem{proposition}[theorem]{Proposition}
\newtheorem{lemma}[theorem]{Lemma}
\newtheorem{corollary}[theorem]{Corollary}
\theoremstyle{definition}
\newtheorem{definition}[theorem]{Definition}

\newcommand{\Z}{\mathbb Z}
\newcommand{\F}{\mathcal F}
\newcommand{\M}{\mathcal M}
\newcommand{\B}{\mathcal B}
\newcommand{\ind}{\mathbf 1}
\numberwithin{equation}{section}
\allowdisplaybreaks[1]
\hypersetup{pdftitle={An Explicit Construction of Large Families of MSTD Sets},pdfauthor={Tahsin Khan}}

\title[Large families of MSTD sets]{An Explicit Construction of Large Families of MSTD Sets}
\author{Tahsin Khan}
\date{}

\begin{document}

\begin{abstract}
A finite set $S\subset\mathbb Z$ is more sums than differences (MSTD) if $|S+S|>|S-S|$. Yufei Zhao's bidirectional-ballot construction gives $\Theta(2^N/N)$ MSTD subsets of $[0,N-1]$, the best previously known size bound for an explicit family, to the author's knowledge. We improve this bound: for every fixed integer $T\geq1$, we construct a family of
\[
 \Theta_T\!\left(\frac{2^N}{\bigl(\log_2^{\circ T}N\bigr)^{2^T}}\right)
\]
MSTD subsets, where $\log_2^{\circ T}$ denotes the $T$-fold iterated logarithm. This family is asymptotically larger than Zhao's for every fixed $T\geq1$. A subfamily of the same order of magnitude can be numbered so that any chosen member can be generated deterministically in time polynomial in $N$, without listing preceding members. As a supplementary comparison, replacing the terminal ballot rule by exhaustive fringe search gives a positive-proportion family with the same polynomial-time generation property.
\end{abstract}
\maketitle

\section{Introduction}\label{sec:intro}

For finite sets $A,B\subset\Z$, write
\[
 A+B=\{a+b:a\in A,\ b\in B\},\qquad
 A-B=\{a-b:a\in A,\ b\in B\}.
\]
A finite set $S\subset\Z$ is an \emph{MSTD set} if $|S+S|>|S-S|$. Throughout, $[a,b]$ denotes the integer interval $\{a,a+1,\ldots,b\}$. We study the size of explicitly described families of MSTD subsets of $[0,N-1]$.

Distinct elements contribute one unordered sum but two signed differences, suggesting that difference sets should often be larger. Nevertheless, MSTD sets exist and admit several types of construction. Nathanson~\cite{Nathanson} developed explicit infinite families by modifying symmetric sets, and Hegarty~\cite{Hegarty} gave further constructions and classified the smallest MSTD sets up to affine transformation. These results demonstrate that sum dominance can be forced by structure; they also motivate the quantitative question of how many distinct sets a construction produces inside a given interval.

The total number of MSTD sets is already known to have order $2^N$. Martin and O'Bryant~\cite{MO} proved that a positive proportion of the subsets of $[0,N-1]$ are MSTD. Zhao~\cite{ZhaoMissing} subsequently proved that this proportion converges to a positive limit and gave an algorithm for approximating it. Thus the problem of constructing large families is distinct from proving that many MSTD sets exist. A useful underlying principle is that extreme sums and differences are controlled by the ends of a set, whereas sums and differences in the middle have many possible representations. Prescribing the fringes and controlling the middle can therefore separate the two cardinalities.

Related probabilistic results clarify the scope of this principle. Hegarty and Miller~\cite{HM} showed that, if each element is included independently with probability $p_N\to0$ and $Np_N\to\infty$, the resulting set is difference-dominated with probability tending to one. Lazarev, Miller, and O'Bryant~\cite{LMO} studied the distribution of missing sums, including bounds for the probability that no possible sum is missing. More recently, He, Daniewski, and Miller~\cite{HDM} considered a Markov-chain model with dependent inclusion decisions and obtained positive limiting probabilities for sum-dominant, difference-dominant, and balanced sets. These counting results concern random families rather than a prescribed combinatorial rule for constructing their members.

For explicit families, Miller, Orosz, and Scheinerman~\cite{MOS} combined fixed fringes with restrictions on gaps in the middle to construct $\Omega(2^N/N^4)$ MSTD subsets of an interval of length $N$. Yufei Zhao~\cite{Zhao} improved this to $\Theta(2^N/N)$ by using bidirectional ballot sequences. Such a sequence has strictly more ones than zeros in every nonempty prefix and suffix. This condition ensures that the associated set fills its possible sum interval; Zhao's enumeration of these sequences then gives the density $\Theta(1/N)$. The survey of Iyer, Lazarev, Miller, and Zhang~\cite{ILMZSurvey} identifies Zhao's construction as the densest explicit family in this line of work. It is the benchmark improved by the present construction.

The literature also includes extensions in other directions. Iyer, Lazarev, Miller, and Zhang~\cite{ILMZGeneral} proved positive-proportion results for generalized comparisons of sums and differences and for sets whose first several iterated sumsets are all sum-dominant. Miller and Peterson~\cite{MP} obtained parametrized families using a geometric description in terms of lattice points in polyhedral cones. Our focus here remains the ordinary comparison $|S+S|>|S-S|$ and the size of a family inside one integer interval.

We construct sets $A\subseteq[0,n]$ satisfying $A+A=[0,2n]$ and then apply Zhao's fixed-fringe conversion to MSTD sets. Instead of requiring a ballot condition on the whole interval, we repeatedly separate it into two short fringes and a middle block. The terminal fringes are bidirectional ballot sets. At each preceding level, the middle subsets satisfy a density condition and a collection of hitting conditions involving the fringes at the next level. These conditions preserve a pairwise sum-filling property while leaving at least one quarter of the possible middle subsets available, uniformly over the deeper choices. Only the terminal fringes incur the ballot-counting loss.

Define
\[
 \ell_0(x)=x,\qquad \ell_{t+1}(x)=\log_2\ell_t(x)
\]
whenever these expressions are defined. Thus $\ell_t(x)$ is the $t$-fold iteration of the base-two logarithm. In all asymptotic statements, $T$ is fixed; constants and thresholds indicated by $O_T$ or $\Theta_T$ may depend on $T$.

\begin{theorem}\label{thm:main}
For every fixed integer $T\geq1$, there is an integer $N_0(T)$ such that, for every $N\geq N_0(T)$, the construction below defines an explicit family $\M_{N,T}$ of MSTD subsets of $[0,N-1]$ with
\[
 |\M_{N,T}|=\Theta_T\!\left(\frac{2^N}{\ell_T(N)^{2^T}}\right).
\]
\end{theorem}

For every fixed $T\geq1$, we have $\ell_T(N)^{2^T}=o(N)$. Consequently, the ratio of the cardinality in Theorem~\ref{thm:main} to the size of Zhao's ballot family has order $N/\ell_T(N)^{2^T}$ and tends to infinity. The improvement concerns the size of a structurally specified family, not the known order of magnitude of the total number of MSTD sets.

Sections~\ref{sec:construction} and~\ref{sec:correctness} give the construction, its cardinality, and the proof that it produces MSTD sets. The appendices give an algorithmic guarantee for this construction and a separate search-based observation; neither is needed for the main proof.

Appendix~\ref{sec:algorithm} asks whether the construction also lets us produce a chosen set efficiently. It numbers a subfamily of the same asymptotic size and, for fixed $T$, constructs the set with a given valid index $i$ in time polynomial in $N$, without generating preceding members. This establishes algorithmic explicitness for a family attaining the size in Theorem~\ref{thm:main}: the combinatorial description supports efficient access to individual sets, not just a count of how many exist.

Appendix~\ref{sec:search-variant} records a supplementary observation: replacing the terminal ballot rule by exhaustive search over short fringe pairs gives a family of size $\Theta(2^N)$, with polynomial-time generation from an index. The proof of the main result does not rely on this variant.

\section{Construction}\label{sec:construction}

\begin{definition}
A set $A\subseteq[0,n]$ is \emph{sum-filled} if $A+A=[0,2n]$.
\end{definition}

Fix $T\geq1$ and define
\begin{equation}\label{eq:lengths}
 k_0=n+1,\qquad
 k_{t+1}=\left\lceil2\log_2 k_t+2t+6\right\rceil
 \quad(0\leq t<T).
\end{equation}

\begin{lemma}\label{lem:scales}
For fixed $T$, as $n\to\infty$,
\[
 k_t=\Theta_T(\ell_t(n))\qquad(0\leq t\leq T).
\]
There is an integer $n_0(T)$ such that, for all $n\geq n_0(T)$,
\begin{equation}\label{eq:separation}
 k_t\geq3k_{t+1}\qquad(0\leq t<T).
\end{equation}
\end{lemma}

\begin{proof}
The estimate for $k_0$ is immediate. At each subsequent level, $k_{t+1}=2\log_2k_t+O_T(1)$. If $k_t=\Theta_T(\ell_t(n))$, then $\log_2k_t=\ell_{t+1}(n)+O_T(1)$, proving the estimate inductively. For each $0\leq t<T$, the ceiling in \eqref{eq:lengths} gives
\[
 k_{t+1}<2\log_2 k_t+2t+7,
 \qquad
 k_t>2^{(k_{t+1}-2t-7)/2}.
\]
Choose $H_t$ such that $2^{(x-2t-7)/2}>3x$ for all $x\geq H_t$. Since $k_{t+1}\to\infty$ as $n\to\infty$, there is an integer $n_t$ such that $k_{t+1}\geq H_t$ whenever $n\geq n_t$. Taking $n_0(T)=\max_{0\leq t<T}n_t$ gives $k_t>3k_{t+1}$ for every $0\leq t<T$ and all $n\geq n_0(T)$.
\end{proof}

Assume henceforth that $n\geq n_0(T)$. Begin with $P_{0,1}=[0,n]$ and $x_{0,1}=0$. At each level $0\leq t\leq T$, there are $2^t$ intervals
\[
 P_{t,a}=[x_{t,a},x_{t,a}+k_t-1],\qquad 1\leq a\leq2^t.
\]
For $0\leq t<T$, partition each into a left fringe, a middle block, and a right fringe:
\begin{align}
 P_{t+1,2a-1}&=[x_{t,a},x_{t,a}+k_{t+1}-1],\nonumber\\
 M_{t,a}&=[x_{t,a}+k_{t+1},x_{t,a}+k_t-k_{t+1}-1],\label{eq:partition}\\
 P_{t+1,2a}&=[x_{t,a}+k_t-k_{t+1},x_{t,a}+k_t-1].\nonumber
\end{align}
Thus
\[
 x_{t+1,2a-1}=x_{t,a},\qquad
 x_{t+1,2a}=x_{t,a}+k_t-k_{t+1}.
\]
Each middle block at level $t$ has length
\begin{equation}\label{eq:m}
 m_t=k_t-2k_{t+1}\geq k_{t+1}.
\end{equation}
The terminal intervals $P_{T,a}$ together with all middle blocks $M_{t,a}$ form a disjoint partition of $[0,n]$. The construction fills these pieces from level $T$ towards level $0$.

\begin{figure}[!htbp]
  \centering
  \includegraphics[width=\linewidth]{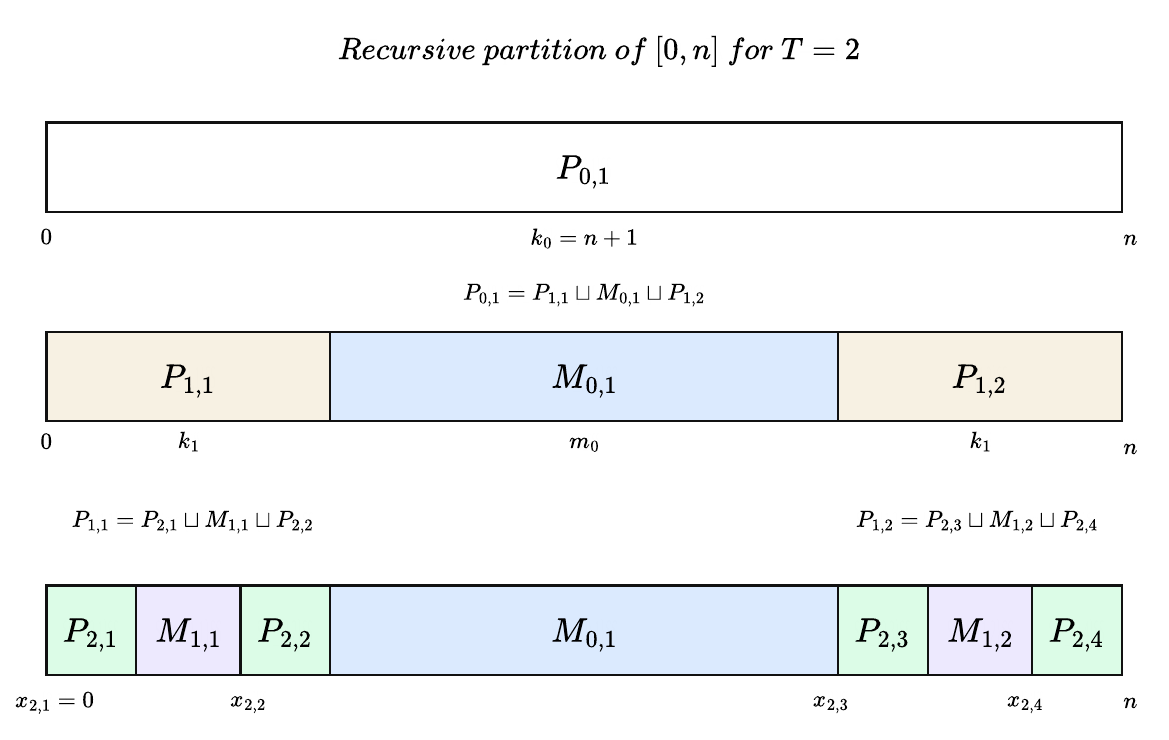}
\caption{Partition of $[0,n]$ for $T=2$. The bottom row comprises the four terminal intervals $P_{2,a}$ and the middle blocks $M_{0,1},M_{1,1},M_{1,2}$. Not to scale. Template prepared with ChatGPT (OpenAI) and edited by the author.}
  \label{fig:recursive-partition}
\end{figure}

Figure~\ref{fig:recursive-partition} illustrates this partition for $T=2$.

Write $P^*_{t,a}\subseteq P_{t,a}$ for the selected subset of a block and
\[
 E_{t,a}=P^*_{t,a}-x_{t,a}\subseteq[0,k_t-1]
\]
for its normalization. Similarly, $M^*_{t,a}\subseteq M_{t,a}$ denotes a selected middle subset. Once the two fringe subsets and the middle subset have been selected, set
\begin{equation}\label{eq:parent}
 P^*_{t,a}=P^*_{t+1,2a-1}\cup M^*_{t,a}\cup P^*_{t+1,2a}.
\end{equation}

The terminal blocks are selected using bidirectional ballot sequences.

\begin{definition}
A binary string is a \emph{bidirectional ballot sequence} if every nonempty prefix and every nonempty suffix has more ones than zeros. A set $E\subseteq[0,h-1]$ is a \emph{bidirectional ballot set} if its indicator string is such a sequence. Let $\B_h$ be the collection of these sets, and let $b_h=|\B_h|$.
\end{definition}

Equivalently, $E\in\B_h$ precisely when, for every $1\leq q\leq h$,
\[
 |E\cap[0,q-1]|>q/2,\qquad
 |E\cap[h-q,h-1]|>q/2.
\]
In particular, each ballot set has more than half of its positions selected. We use Zhao's enumeration result:
\begin{equation}\label{eq:bbs-count}
 b_h=\Theta(2^h/h)\qquad(h\to\infty)
 \quad\text{\cite[Proposition 3.2]{Zhao}}.
\end{equation}
At the terminal level, choose independently
\[
 E_{T,a}\in\B_{k_T},\qquad
 P^*_{T,a}=x_{T,a}+E_{T,a}
 \quad(1\leq a\leq2^T).
\]
There are exactly $b_{k_T}^{\,2^T}$ terminal configurations.

For $t=T-1,T-2,\ldots,0$, suppose the selected blocks at level $t+1$ have been fixed. For each $1\leq a\leq2^t$, call a middle subset $M^*_{t,a}$ \emph{admissible} if it satisfies the density condition
\begin{equation}\tag{D}\label{eq:D}
 |M^*_{t,a}|\geq\left\lceil\frac{m_t}{2}\right\rceil
\end{equation}
\begin{samepage}
and the hitting conditions
\begin{equation}\tag{H}\label{eq:H}
 M^*_{t,a}\cap\bigl(x_{t,a}+r-E_{t+1,b}\bigr)\neq\varnothing
\end{equation}
for all
\begin{equation}\label{eq:H-range}
 1\leq b\leq2^{t+1},\qquad
 2k_{t+1}-1\leq r\leq k_t-k_{t+1}-1.
\end{equation}
\end{samepage}
For every $u\in E_{t+1,b}\subseteq[0,k_{t+1}-1]$ and $r$ in \eqref{eq:H-range},
\[
 x_{t,a}+k_{t+1}
 \leq x_{t,a}+r-u
 \leq x_{t,a}+k_t-k_{t+1}-1.
\]
Thus $x_{t,a}+r-E_{t+1,b}$ is a subset of $M_{t,a}$, so every position required by \eqref{eq:H} lies within the middle interval.

For each level $0\leq t<T$, set
\begin{equation}\label{eq:q}
 q_t=2^{m_t-2}.
\end{equation}
Lemma~\ref{lem:quarter} in Section~\ref{sec:correctness} shows that every configuration produced at deeper levels leaves at least $q_t$ admissible choices for each middle block, and that every selected block contains more than half of its positions.

\begin{definition}\label{def:family}
The family $\F_{n,T}$ consists of all sets obtained by choosing the normalized terminal sets from $\B_{k_T}$ and then, successively for $t=T-1,\ldots,0$, choosing any admissible middle subset in each block. Form each $P^*_{t,a}$ by \eqref{eq:parent}; the final set is $A=P^*_{0,1}$.
\end{definition}

Let
\begin{equation}\label{eq:QJ}
 J=\sum_{t=0}^{T-1}2^t=2^T-1,\qquad
 Q=\sum_{t=0}^{T-1}2^tm_t=n+1-2^Tk_T.
\end{equation}
Here $J$ is the number of middle blocks. Since the middle blocks and the $2^T$ terminal intervals partition $[0,n]$, their total length is $Q=n+1-2^Tk_T$.

\begin{proposition}\label{prop:count}
For fixed $T$, the cardinality of $\F_{n,T}$ satisfies
\begin{equation}\label{eq:asymptotic}
 |\F_{n,T}|=\Theta_T\!\left(\frac{2^{n+1}}{k_T^{2^T}}\right)
 =\Theta_T\!\left(\frac{2^{n+1}}{\ell_T(n)^{2^T}}\right).
\end{equation}
\end{proposition}

\begin{proof}
There are $b_{k_T}^{\,2^T}$ terminal configurations. For any fixed configuration of the deeper levels, each of the $2^t$ middle blocks at level $t$ has at least $q_t$ admissible choices by Lemma~\ref{lem:quarter}. Their conditions do not involve one another, so there are at least $q_t^{2^t}$ extensions at that level. The actual numbers may depend on the deeper choices, but these lower bounds do not. They can therefore be multiplied over the levels.

There is no overcounting. The terminal and middle intervals partition $[0,n]$, so the final set determines the subset selected in each piece and hence the entire sequence of choices. Therefore,
\begin{align}
 |\F_{n,T}|&\geq V_{n,T}:=
 b_{k_T}^{\,2^T}\prod_{t=0}^{T-1}q_t^{\,2^t}\nonumber\\
 &=b_{k_T}^{\,2^T}\,2^{\,n+1-2^Tk_T-2(2^T-1)}.
 \label{eq:count-bound}
\end{align}
For the final equality, \eqref{eq:q} and \eqref{eq:QJ} give
\[
 \prod_{t=0}^{T-1}q_t^{\,2^t}
 =2^{\sum_{t=0}^{T-1}2^t(m_t-2)}
 =2^{Q-2J}
 =2^{n+1-2^Tk_T-2(2^T-1)}.
\]
Applying \eqref{eq:bbs-count}, with $T$ fixed, yields
\[
 V_{n,T}
 =\Theta_T\!\left[
   \left(\frac{2^{k_T}}{k_T}\right)^{2^T}
   2^{n+1-2^Tk_T-2(2^T-1)}
 \right]
 =\Theta_T\!\left(\frac{2^{n+1}}{k_T^{2^T}}\right).
\]
The factor $2^{-2(2^T-1)}$ depends only on $T$. Since $|\F_{n,T}|\geq V_{n,T}$, this proves the required lower bound.

For the upper bound, each middle block of length $m_t$ admits at most $2^{m_t}$ choices, regardless of the selected deeper blocks. Hence, using \eqref{eq:QJ},
\begin{align*}
 |\F_{n,T}|
 &\leq b_{k_T}^{\,2^T}\prod_{t=0}^{T-1}2^{\,2^tm_t}
 =b_{k_T}^{\,2^T}\,2^{\,n+1-2^Tk_T}\\
 &=\Theta_T\!\left(\frac{2^{n+1}}{k_T^{2^T}}\right),
\end{align*}
where the final equality again follows from \eqref{eq:bbs-count}. Together with the lower bound, this proves the first equality in \eqref{eq:asymptotic}. Lemma~\ref{lem:scales} gives the iterated-logarithm form.
\end{proof}

\section{Sum-filling property}\label{sec:correctness}

We prove that every set obtained from the construction is sum-filled.

\begin{lemma}\label{lem:ballot-sums}
If $E,F\in\B_h$, then $E+F=[0,2h-2]$.
\end{lemma}

\begin{proof}
For $0\leq s\leq h-1$, the sets
\[
 E\cap[0,s],\qquad s-(F\cap[0,s])
\]
each occupy more than half of $[0,s]$. They intersect, so $s=e+f$ for some $e\in E$ and $f\in F$. The reflected sets $(h-1)-E$ and $(h-1)-F$ are also ballot sets. Applying the same argument to them supplies the remaining sums $h,\ldots,2h-2$.
\end{proof}

\begin{proposition}\label{prop:pairwise}
For every $A\in\F_{n,T}$, all $0\leq t\leq T$ and $1\leq a,b\leq2^t$ satisfy
\begin{equation}\label{eq:pairwise}
 P^*_{t,a}+P^*_{t,b}
 =[x_{t,a}+x_{t,b},\ x_{t,a}+x_{t,b}+2k_t-2].
\end{equation}
\end{proposition}

\begin{proof}
We proceed from $t=T$ down to $t=0$. At the terminal level, the claim follows from Lemma~\ref{lem:ballot-sums} by translation.

Suppose it holds at level $t+1$. Fix $a,b$ at level $t$ and put
\[
 K=k_t,\qquad h=k_{t+1},\qquad \beta=x_{t,a}+x_{t,b}.
\]
The target interval $[\beta,\beta+2K-2]$ is the disjoint union of the following five consecutive intervals:
\begin{align*}
 I_1&=[\beta,\beta+2h-2],\\
 I_2&=[\beta+2h-1,\beta+K-h-1],\\
 I_3&=[\beta+K-h,\beta+K+h-2],\\
 I_4&=[\beta+K+h-1,\beta+2K-2h-1],\\
 I_5&=[\beta+2K-2h,\beta+2K-2].
\end{align*}
The condition $K\geq3h$ ensures these intervals are well-defined. Figure~\ref{fig:five-intervals} illustrates the case $t=0$ and $a=b=1$.

\begin{figure}[!htbp]
  \centering
  \includegraphics[width=\linewidth]{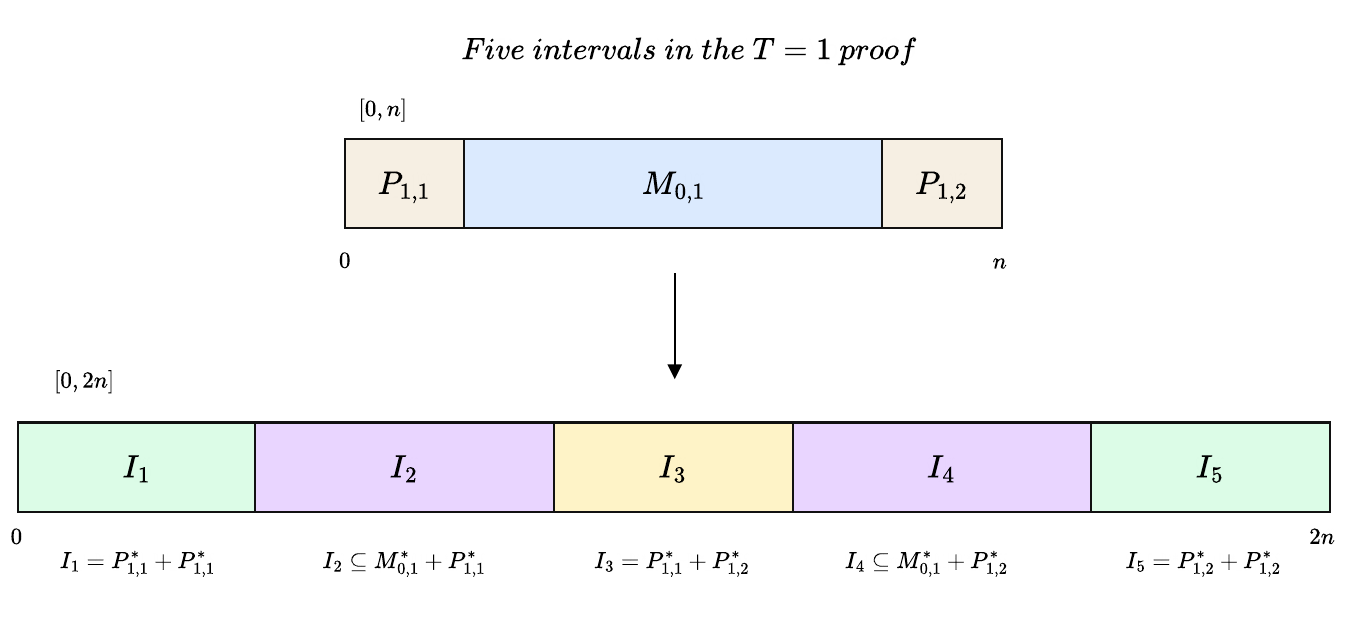}
\caption{Top-level sum-filling: fringe sums give $I_1,I_3,I_5$; hitting conditions cover $I_2,I_4$. Stars denote selected subsets. Not to scale. Template prepared with ChatGPT (OpenAI) and edited by the author.}
  \label{fig:five-intervals}
\end{figure}

The induction hypothesis supplies three of the intervals:
\begin{align*}
 P^*_{t+1,2a-1}+P^*_{t+1,2b-1}&=I_1,\\
 P^*_{t+1,2a-1}+P^*_{t+1,2b}&=I_3,\\
 P^*_{t+1,2a}+P^*_{t+1,2b}&=I_5.
\end{align*}
All these sets are contained in $P^*_{t,a}+P^*_{t,b}$.

For $s\in I_2$, write $s=\beta+r$, where $2h-1\leq r\leq K-h-1$. Condition \eqref{eq:H}, with index $2b-1$, supplies $u\in E_{t+1,2b-1}$ with
\[
 q=x_{t,a}+r-u\in M^*_{t,a}.
\]
The left fringe of $P_{t,b}$ begins at $x_{t,b}$, so $p=x_{t,b}+u\in P^*_{t+1,2b-1}$. Hence $p+q=s$.

For $s\in I_4$, write $s=\beta+K-h+r$ with the same range of $r$. Using index $2b$ in \eqref{eq:H} gives $u\in E_{t+1,2b}$ such that $q=x_{t,a}+r-u\in M^*_{t,a}$. The right fringe of $P_{t,b}$ begins at $x_{t,b}+K-h$, so
\[
 p=x_{t,b}+K-h+u\in P^*_{t+1,2b},\qquad p+q=s.
\]
Thus the two remaining gaps are also covered. The reverse containment in \eqref{eq:pairwise} follows from the ambient intervals of the selected blocks, completing the induction.
\end{proof}

\begin{corollary}\label{cor:sumfilled}
Every $A\in\F_{n,T}$ satisfies $A+A=[0,2n]$.
\end{corollary}

\begin{proof}
Apply Proposition~\ref{prop:pairwise} with $t=0$ and $a=b=1$.
\end{proof}

We use Zhao's fixed fringes~\cite{Zhao} to convert sum-filled sets into MSTD sets.

\begin{lemma}\label{lem:fringes}
Let $A\subseteq[0,n]$ be sum-filled and set $N=n+23$. Define
\begin{align*}
 L_0&=\{0,2,3,7,8,9,10\},\\
 R_N&=\{N-11,N-10,N-9,N-8,N-6,N-3,N-2,N-1\}.
\end{align*}
Then
\[
 \Psi_n(A)=L_0\cup(A+11)\cup R_N
\]
is an MSTD subset of $[0,N-1]$. For fixed $n$, the map $\Psi_n$ is injective.
\end{lemma}

\begin{proof}
Put $S=\Psi_n(A)$. Since the extreme sums $0$ and $2n$ belong to $A+A$, we have $0,n\in A$. The three parts of $S$ lie in the disjoint intervals $[0,10]$, $[11,N-12]$, and $[N-11,N-1]$. Direct addition gives
\[
 L_0+L_0=[0,20]\setminus\{1\},\qquad
 R_N+R_N=[2N-22,2N-2].
\]
Also,
\[
 (A+11)+(A+11)=[22,2N-24],
\]
and the two gaps between these intervals are covered by
\[
 21=10+11,\qquad 2N-23=(N-11)+(N-12).
\]
The sum $1$ is absent, because $S$ is nonnegative and $1\notin S$. Consequently,
\begin{equation}\label{eq:sum-fringe}
 S+S=[0,2N-2]\setminus\{1\},\qquad |S+S|=2N-2.
\end{equation}

If $a-b=N-7$ with $a,b\in S$, then $0\leq b\leq6$. The only possible values of $b$ are $0,2,3$, which would require $a\in\{N-7,N-5,N-4\}$. None belongs to $S$. Thus both $N-7$ and $-(N-7)$ are absent from $S-S$, giving
\[
 |S-S|\leq2N-3<2N-2=|S+S|.
\]
Finally, $A=(S\cap[11,N-12])-11$, proving injectivity.
\end{proof}

For $N\geq n_0(T)+23$, define
\begin{equation}\label{eq:M-family}
 \M_{N,T}=\{\Psi_{N-23}(A):A\in\F_{N-23,T}\}.
\end{equation}
By Lemma~\ref{lem:fringes} and Corollary~\ref{cor:sumfilled}, its members are MSTD and the conversion is injective. Therefore
\begin{align}
 |\M_{N,T}|&=|\F_{N-23,T}|\nonumber\\
 &=\Theta_T\!\left(\frac{2^{N-22}}{\ell_T(N-23)^{2^T}}\right)\nonumber\\
 &=\Theta_T\!\left(\frac{2^N}{\ell_T(N)^{2^T}}\right).
 \label{eq:M-count}
\end{align}
where the last step uses $\ell_T(N-23)/\ell_T(N)\to1$ for fixed $T$.

The counting argument in Section~\ref{sec:construction} rests on the following availability statement.

\begin{lemma}\label{lem:quarter}
Suppose every selected block at level $t+1$ contains more than half of its positions. For every middle block at level $t$, at least $2^{m_t-2}$ of its subsets satisfy both \eqref{eq:D} and \eqref{eq:H}. Every set $P^*_{t,a}$ formed using such a middle subset also contains more than half of the positions in $P_{t,a}$.
\end{lemma}

\begin{proof}
Write $K=k_t$, $h=k_{t+1}$, and $m=m_t$. Fix a middle block. For one pair $(b,r)$, failure of \eqref{eq:H} means omitting the $|E_{t+1,b}|>h/2$ elements of its translate. The number of failing subsets is
\[
 2^{m-|E_{t+1,b}|}<2^{m-h/2}.
\]
There are $2^{t+1}$ choices of $b$ and $K-3h+1\leq K$ choices of $r$. By the union bound, the number of subsets failing at least one hitting condition is at most
\begin{equation}\label{eq:bad}
 2^{t+1}K\,2^{m-h/2}\leq2^{m-2},
\end{equation}
where the last inequality uses $h\geq2\log_2K+2t+6$.

At least $2^{m-1}$ subsets have cardinality at least $\lceil m/2\rceil$. Subtracting the upper bound in \eqref{eq:bad} leaves at least
\[
 2^{m-1}-2^{m-2}=2^{m-2}
\]
subsets satisfying both conditions.

For any such choice, the disjoint union \eqref{eq:parent} gives
\[
 |P^*_{t,a}|>
 \frac h2+\frac{K-2h}{2}+\frac h2=\frac K2.
\]
\end{proof}

Starting with the terminal ballot sets, apply Lemma~\ref{lem:quarter} successively for $t=T-1,T-2,\ldots,0$. Every deeper configuration then has the required number of admissible extensions, and every selected block satisfies
\begin{equation}\label{eq:strict-density}
 |P^*_{t,a}|>k_t/2.
\end{equation}
Together with \eqref{eq:M-count}, this completes the proof of Theorem~\ref{thm:main}.

\appendix
\section{Polynomial-time generation of individual sets}\label{sec:algorithm}

This appendix supplies the algorithmic guarantee accompanying the main construction: a subfamily of the same asymptotic size can be numbered so that any chosen set is recovered from its index in polynomial time, without generating earlier sets. The proof below gives the explicit bound $O_T(n^5)$ on bit operations, meaning elementary operations on binary digits.

For each configuration at level $t+1$, retain only the first $q_t$ admissible middle subsets in lexicographic order of their left-to-right indicator strings, with $0<1$. This defines a subfamily $\widehat{\F}_{n,T}\subseteq\F_{n,T}$. Order the terminal ballot sets in the same way. All indices in this appendix start at zero.

\begin{theorem}\label{thm:indexing}
For every fixed $T\geq1$ and $n\geq n_0(T)$, the subfamily $\widehat{\F}_{n,T}$ has cardinality
\[
 |\widehat{\F}_{n,T}|=V_{n,T}
 =\Theta_T\!\left(\frac{2^{n+1}}{\ell_T(n)^{2^T}}\right)
\]
and admits a bijective indexing by $0\leq i<V_{n,T}$. Its $i$-th member can be computed deterministically in $O_T(n^5)$ bit operations, using dynamic programming.
\end{theorem}

\begin{proof}
Lemma~\ref{lem:quarter} ensures that the retained choices exist. There are $b_{k_T}^{\,2^T}$ terminal configurations and now exactly $q_t$ choices in each middle block at level $t$. The multiplication and no-overcounting argument of Proposition~\ref{prop:count} gives $|\widehat{\F}_{n,T}|=V_{n,T}$, with the stated asymptotic size.

First enumerate $\B_{k_T}$ and compute $b_{k_T}$. Define
\begin{equation}\label{eq:rank-counts}
 f_T=b_{k_T},\qquad f_t=q_t f_{t+1}^2\quad(0\leq t<T).
\end{equation}
Here $f_t$ counts assignments of local ranks to a level-$t$ block and all the blocks within it. In particular,
\[
 f_0=b_{k_T}^{\,2^T}\prod_{t=0}^{T-1}q_t^{\,2^t}=V_{n,T}.
\]
Starting with index $i$ at $P_{0,1}$, split an index $0\leq z<f_t$ at any nonterminal block into
\begin{equation}\label{eq:recursive-ranks}
 \begin{aligned}
 j&=z\bmod q_t,& v&=\lfloor z/q_t\rfloor,\\
 z_L&=v\bmod f_{t+1},& z_R&=\lfloor v/f_{t+1}\rfloor.
 \end{aligned}
\end{equation}
Then $0\leq j<q_t$ and $0\leq z_L,z_R<f_{t+1}$. Assign $j$ to its middle block and continue recursively with $z_L,z_R$ on the two fringes. At level $T$, the remaining indices select the terminal ballot sets.

This recursion assigns ranks only. The sets are constructed afterwards, from level $T$ upwards, because \eqref{eq:H} involves every block at the next level, including those outside the current fringe. Thus the two fringes are not constructed independently. Once all selected blocks at level $t+1$ are known, each middle subset at level $t$ is chosen using its assigned rank as follows.

Write $K=k_t$, $h=k_{t+1}$, $m=m_t$, and $d=\lceil m/2\rceil$. Let $E_1,\ldots,E_g$, where $g=2^{t+1}$, be all the normalized sets at level $t+1$. Encode a middle subset by $y=(y_0,\ldots,y_{m-1})$, where $y_p=1$ means that $x_{t,a}+h+p$ is selected. Conditions \eqref{eq:H} and \eqref{eq:D} become
\begin{equation}\label{eq:local-conditions}
 \sum_{u\in E_b}y_{p-u}\geq1
 \quad\left(1\leq b\leq g,\ h-1\leq p\leq m-1\right),
 \qquad \sum_{p=0}^{m-1}y_p\geq d.
\end{equation}
Each hitting condition depends on a window of $h$ consecutive bits.

For a prefix of length $p$ whose completed windows pass these tests, let $s$ be its last $\min(p,h-1)$ bits and $c$ its total number of ones. Let $C(p,s,c)$ count its valid completions. If a bit $\varepsilon\in\{0,1\}$ is appended, let $s_\varepsilon$ be the last at most $h-1$ bits of $s\varepsilon$. Call the bit permitted if the new length-$h$ window, when there is one, satisfies all the hitting tests in \eqref{eq:local-conditions}. Then
\begin{equation}\label{eq:completion-dp}
 \begin{aligned}
 C(m,s,c)&=\ind_{\{c\geq d\}},\\
 C(p,s,c)&=\sum_{\varepsilon\text{ permitted}}
 C(p+1,s_\varepsilon,c+\varepsilon)\qquad(p<m).
 \end{aligned}
\end{equation}
The two next-bit classes are disjoint, and future conditions depend only on the updated state. Therefore the recurrence, evaluated in decreasing order of $p$, counts exactly the completions satisfying both density and hitting conditions.

To obtain the middle subset of rank $j$, start with the empty prefix. Let $C_0=C(p+1,s_0,c)$ if appending $0$ is permitted, and $C_0=0$ otherwise. If $j<C_0$, append $0$; otherwise append $1$ and replace $j$ by $j-C_0$. Update $p,s,c$ and repeat for all $m$ positions. This selects the admissible word of the assigned rank, since the $C_0$ completions beginning with $0$ precede those beginning with $1$. Lemma~\ref{lem:quarter} guarantees that every assigned rank below $q_t$ exists. The normalized tests are identical for all middle blocks at the same level, so one completion table serves them all. Form each $P^*_{t,a}$ using \eqref{eq:parent}, and continue upwards until $A=P^*_{0,1}$ is obtained.

The rank splitting is invertible, since
\[
 z=j+q_t(z_L+f_{t+1}z_R).
\]
Every rank assignment produces a set in $\widehat{\F}_{n,T}$. Conversely, a set in this subfamily determines its terminal ballot ranks and, level by level, the rank of each middle subset in the appropriate list. The disjoint partition prevents any of these choices from being lost. Hence the indexing is bijective.

We charge $O(b)$ bit operations for addition, comparison, copying, or subtraction of $b$-bit integers. At level $t$,
\[
 2^h<2^{2t+7}K^2=O_T(K^2).
\]
All permitted windows can be precomputed in $O_T(h2^h)$ bit operations. The completion table has $O(m^2 2^h)=O_T(K^4)$ states, each with two possible transitions and counts of at most $m+1=O(K)$ bits. Thus it costs $O_T(K^5)$ bit operations. Decoding all middle subsets after their tables are available costs $O(\sum_{t<T}2^t m_t^2)=O(n^2)$ bit operations.

Enumerating the terminal ballot sets costs $O(k_T\log k_T\,2^{k_T})$ bit operations, by checking prefixes and suffixes of every length-$k_T$ word. Since $2^{k_T}=O_T(k_{T-1}^2)$ and $k_T=O_T(\log n)$, this is at most $O_T(n^2\log n\log\log n)$. The numbers in the rank recursion have $O(n)$ bits, and there are $O_T(1)$ blocks, so computing the $f_t$ and splitting the ranks costs $O_T(n^2)$ bit operations using elementary binary arithmetic. The lengths in \eqref{eq:lengths} can also be computed exactly, using $\lceil2\log_2 k_t\rceil=\lceil\log_2(k_t^2)\rceil$. The total running time is therefore
\[
 O_T\!\left(\sum_{t=0}^{T-1}k_t^5\right)=O_T(n^5).\qedhere
\]
\end{proof}

For $N=n+23$, applying $\Psi_n$ gives the corresponding indexed member of the subfamily
\[
 \widehat{\M}_{N,T}
 =\{\Psi_{N-23}(A):A\in\widehat{\F}_{N-23,T}\}
 \subseteq\M_{N,T}.
\]
The map is injective, so this subfamily has size $V_{N-23,T}=\Theta_T(2^N/\ell_T(N)^{2^T})$. The conversion adds only $O(N)$ bit operations, giving $O_T(N^5)$ bit operations per indexed member.

\section{A larger family using fringe search}\label{sec:search-variant}

This appendix records a supplementary search-based variant. It uses Proposition~\ref{prop:pairwise} and Lemma~\ref{lem:quarter}, but the main proof does not rely on this variant. Instead of starting from bidirectional ballot sets, we find compatible fringe pairs by exhaustive search. Already for $T=1$, this yields $\Theta(2^n)$ sum-filled sets, with $O(n^5)$ bit operations to generate a chosen member.

Put
\[
 h=k_1=\left\lceil2\log_2(n+1)+6\right\rceil,\qquad
 m=m_0=n+1-2h,\qquad q=2^{m-2},
\]
and take $n$ sufficiently large that $n+1\geq3h$. Enumerate all ordered pairs $(E,F)$ of subsets of $[0,h-1]$, retaining those satisfying
\begin{equation}\label{eq:fringe-compatible}
 |E|,|F|>h/2,\qquad E+E=E+F=F+F=[0,2h-2].
\end{equation}
Let $c_h$ be the number of retained pairs.

The fixed-fringe argument of Martin and O'Bryant~\cite[Proposition~8]{MO}, used for sum-filled sets in~\cite[Section~4]{LMO}, adapts to give $c_h=\Theta(2^{2h})$. We briefly check the additional conditions. For a fixed positive integer $L$ and $h>2L$, require both sets to contain $[0,L-1]\cup[h-L,h-1]$, and choose all remaining positions independently with probability $1/2$. Both density inequalities then hold with probability at least $1/4$. The forced positions fill the extreme sums in all three sumsets. For any remaining sum $s$, put $d_s=\min(s,2h-2-s)$. In each sumset, including $E+F$, there are $\lfloor(d_s+1)/2\rfloor$ disjoint pairs of distinct positions giving independent representations, each present with probability at least $1/4$. A union bound therefore gives
\[
 \Pr(\text{at least one of the three sumsets is incomplete})
 \leq 6\sum_{d=2L-1}^{\infty}(3/4)^{\lfloor(d+1)/2\rfloor}
 =48(3/4)^L.
\]
These probabilities are conditional on the forced positions. Choose $L$ so that the last quantity is less than $1/8$. At least one eighth of the $2^{2h-4L}$ conditioned pairs satisfy all of \eqref{eq:fringe-compatible}. This proves the lower bound for $c_h$; the upper bound $c_h\leq2^{2h}$ is immediate.

For each retained pair, keep the first $q$ middle subsets in lexicographic order, with $M\subseteq[h,n-h]$ satisfying \eqref{eq:D} and \eqref{eq:H} at $t=0$. Lemma~\ref{lem:quarter} guarantees at least $q$ choices, and the single-level case of Proposition~\ref{prop:pairwise} shows that
\[
 A=E\cup M\cup(n+1-h+F)
\]
is sum-filled. The three pieces determine their choices uniquely, so the resulting family $\mathcal G_n$ satisfies
\[
 |\mathcal G_n|=c_hq=\Theta(2^{2h})\,2^{n+1-2h-2}
 =\Theta(2^n).
\]
The fixed-fringe map $\Psi_n$ gives $\Theta(2^N)$ MSTD subsets of $[0,N-1]$, where $N=n+23$.

Order the retained fringe pairs lexicographically by their concatenated indicator strings. Given $0\leq i<c_hq$, write
\[
 j=i\bmod q,\qquad i'=\frac{i-j}{q}.
\]
Take the zero-based $i'$-th retained pair and use the algorithm of Appendix~\ref{sec:algorithm} to obtain its $j$-th admissible middle subset. This gives a bijective indexing of $\mathcal G_n$, with the same division of the index into fringe and middle ranks as before.

There are $2^{2h}=\Theta(n^4)$ candidate fringe pairs. Directly testing a pair costs $O(h^2\log h)$ bit operations, so constructing the ordered list costs $O(n^4(\log n)^2\log\log n)$. The middle-block algorithm costs $O(n^5)$ bit operations and dominates this search. Thus one indexed member is produced in $O(n^5)$ bit operations; this is not a bound on the time to list the entire family.

\section*{Acknowledgments}

The author thanks Steven J. Miller for his comments on an early version of the work and his suggestion to pursue publication, and Evan Chen for helpful suggestions on the presentation of the manuscript.

\section*{Declaration of AI assistance}

ChatGPT (OpenAI) suggested using different fringe sets rather than a fixed fringe set and assisted with developing the algorithm in Appendix~\ref{sec:algorithm}, manuscript preparation, and draw.io diagram templates. Aristotle AI was used for automated feedback on the arguments. The author takes responsibility for the final mathematical content and manuscript.


\begin{thebibliography}{99}

\bibitem{HDM}
F.~He, K.~Daniewski, and S.~J.~Miller, \emph{Generating and generalizing MSTD sets through Markov processes}, preprint (2026), \href{https://arxiv.org/abs/2610.08394}{arXiv:2610.08394}.

\bibitem{Hegarty}
P.~V.~Hegarty, \emph{Some explicit constructions of sets with more sums than differences}, Acta Arithmetica \textbf{130} (2007), no.~1, 61--77. \href{https://doi.org/10.4064/aa130-1-4}{doi:10.4064/aa130-1-4}.

\bibitem{HM}
P.~Hegarty and S.~J.~Miller, \emph{When almost all sets are difference dominated}, Random Structures \& Algorithms \textbf{35} (2009), no.~1, 118--136. \href{https://doi.org/10.1002/rsa.20268}{doi:10.1002/rsa.20268}.

\bibitem{ILMZGeneral}
G.~Iyer, O.~Lazarev, S.~J.~Miller, and L.~Zhang, \emph{Generalized more sums than differences sets}, Journal of Number Theory \textbf{132} (2012), no.~5, 1054--1073. \href{https://arxiv.org/abs/1108.4500}{arXiv:1108.4500}.

\bibitem{ILMZSurvey}
G.~Iyer, O.~Lazarev, S.~J.~Miller, and L.~Zhang, \emph{Finding and counting MSTD sets}, in \emph{Combinatorial and Additive Number Theory: CANT 2011 and 2012} (M.~B.~Nathanson, ed.), Springer Proceedings in Mathematics \& Statistics, vol.~101, Springer, New York, 2014, pp.~79--98. \href{https://doi.org/10.1007/978-1-4939-1601-6_7}{doi:10.1007/978-1-4939-1601-6\_7}.

\bibitem{LMO}
O.~Lazarev, S.~J.~Miller, and K.~O'Bryant, \emph{Distribution of missing sums in sumsets}, Experimental Mathematics \textbf{22} (2013), no.~2, 132--156. \href{https://arxiv.org/abs/1109.4700}{arXiv:1109.4700}.

\bibitem{MO}
G.~Martin and K.~O'Bryant, \emph{Many sets have more sums than differences}, in \emph{Additive Combinatorics}, CRM Proceedings \& Lecture Notes, vol.~43, American Mathematical Society, Providence, RI, 2007, pp.~287--305. \href{https://arxiv.org/abs/math/0608131}{arXiv:math/0608131}.

\bibitem{MOS}
S.~J.~Miller, B.~Orosz, and D.~Scheinerman, \emph{Explicit constructions of infinite families of MSTD sets}, Journal of Number Theory \textbf{130} (2010), no.~5, 1221--1233. \href{https://doi.org/10.1016/j.jnt.2009.09.003}{doi:10.1016/j.jnt.2009.09.003}.

\bibitem{MP}
S.~J.~Miller and C.~Peterson, \emph{A geometric perspective on the MSTD question}, Discrete \& Computational Geometry \textbf{62} (2019), no.~4, 832--855. \href{https://doi.org/10.1007/s00454-019-00109-7}{doi:10.1007/s00454-019-00109-7}.

\bibitem{Nathanson}
M.~B.~Nathanson, \emph{Sets with more sums than differences}, Integers \textbf{7} (2007), Paper A5, 24~pp. \href{https://arxiv.org/abs/math/0608148}{arXiv:math/0608148}.

\bibitem{Zhao}
Y.~Zhao, \emph{Constructing MSTD sets using bidirectional ballot sequences}, Journal of Number Theory \textbf{130} (2010), no.~5, 1212--1220. \href{https://doi.org/10.1016/j.jnt.2009.11.005}{doi:10.1016/j.jnt.2009.11.005}.

\bibitem{ZhaoMissing}
Y.~Zhao, \emph{Sets characterized by missing sums and differences}, Journal of Number Theory \textbf{131} (2011), no.~11, 2107--2134. \href{https://doi.org/10.1016/j.jnt.2011.05.003}{doi:10.1016/j.jnt.2011.05.003}.

\end{thebibliography}
\end{document}